\documentclass[12pt]{article}

\catcode`\@=11

\thicklines
\newskip\Einheit \Einheit=.6cm
\newcount\xcoord \newcount\ycoord
\newdimen\xdim \newdimen\ydim \newdimen\PfadD@cke \newdimen\Pfadd@cke
\def\PfadDicke#1{\PfadD@cke#1 \divide\PfadD@cke by2 
\Pfadd@cke\PfadD@cke \multiply\PfadD@cke by2}
\long\def\LOOP#1\REPEAT{\def\BODY{#1}\ITERATE}
\def\ITERATE{\BODY \let\next\ITERATE \else\let\next\relax\fi \next}
\let\REPEAT=\fi
\def\Punkt{\hbox{\raise-2pt\hbox to0pt{\hss\scriptsize$\bullet$\hss}}}

\def\DuennPunkt(#1,#2){\unskip
  \raise#2 \Einheit\hbox to0pt{\hskip#1 \Einheit
          \raise-1.5pt\hbox to0pt{\hss\tiny$\bullet$\hss}\hss}}
		  
\def\NormalPunkt(#1,#2){\unskip
  \raise#2 \Einheit\hbox to0pt{\hskip#1 \Einheit
          \raise-3pt\hbox to0pt{\hss\large$\bullet$\hss}\hss}}
\def\DickPunkt(#1,#2){\unskip
  \raise#2 \Einheit\hbox to0pt{\hskip#1 \Einheit
          \raise-4pt\hbox to0pt{\hss\Large$\bullet$\hss}\hss}}
\def\Kreis(#1,#2){\unskip
  \raise#2 \Einheit\hbox to0pt{\hskip#1 \Einheit
          \raise-4pt\hbox to0pt{\hss\Large$\circ$\hss}\hss}}
\def\Diagonale(#1,#2)#3{\unskip\leavevmode
  \xcoord#1\relax \ycoord#2\relax
      \raise\ycoord \Einheit\hbox to0pt{\hskip\xcoord \Einheit
         \unitlength\Einheit
         \line(1,1){#3}\hss}}
\def\AntiDiagonale(#1,#2)#3{\unskip\leavevmode
  \xcoord#1\relax \ycoord#2\relax \advance\xcoord by -0.05\relax
      \raise\ycoord \Einheit\hbox to0pt{\hskip\xcoord \Einheit
         \unitlength\Einheit
         \line(1,-1){#3}\hss}}
\def\Pfad(#1,#2),#3\endPfad{\unskip\leavevmode
  \xcoord#1 \ycoord#2 \thicklines\ZeichnePfad#3\endPfad\thinlines}
\def\ZeichnePfad#1{\ifx#1\endPfad\let\next\relax
  \else\let\next\ZeichnePfad
    \ifnum#1=1
      \raise\ycoord \Einheit\hbox to0pt{\hskip\xcoord \Einheit
         \vrule height\Pfadd@cke width1 \Einheit depth\Pfadd@cke\hss}%
      \advance\xcoord by 1
    \else\ifnum#1=2
      \raise\ycoord \Einheit\hbox to0pt{\hskip\xcoord \Einheit
        \hbox{\hskip-\PfadD@cke\vrule height1 \Einheit 
width\PfadD@cke depth0pt}\hss}%
      \advance\ycoord by 1
    \else\ifnum#1=3
      \raise\ycoord \Einheit\hbox to0pt{\hskip\xcoord \Einheit
         \unitlength\Einheit
         \line(1,1){1}\hss}
      \advance\xcoord by 1
      \advance\ycoord by 1
    \else\ifnum#1=4
      \raise\ycoord \Einheit\hbox to0pt{\hskip\xcoord \Einheit
         \unitlength\Einheit
         \line(1,-1){1}\hss}
      \advance\xcoord by 1
      \advance\ycoord by -1
    \fi\fi\fi\fi
  \fi\next}
\def\hSSchritt{\leavevmode\raise-.4pt\hbox 
to0pt{\hss.\hss}\hskip.2\Einheit
  \raise-.4pt\hbox to0pt{\hss.\hss}\hskip.2\Einheit
  \raise-.4pt\hbox to0pt{\hss.\hss}\hskip.2\Einheit
  \raise-.4pt\hbox to0pt{\hss.\hss}\hskip.2\Einheit
  \raise-.4pt\hbox to0pt{\hss.\hss}\hskip.2\Einheit}
\def\vSSchritt{\vbox{\baselineskip.2\Einheit\lineskiplimit0pt
\hbox{.}\hbox{.}\hbox{.}\hbox{.}\hbox{.}}}
\def\DSSchritt{\leavevmode\raise-.4pt\hbox to0pt{%
  \hbox to0pt{\hss.\hss}\hskip.2\Einheit
  \raise.2\Einheit\hbox to0pt{\hss.\hss}\hskip.2\Einheit
  \raise.4\Einheit\hbox to0pt{\hss.\hss}\hskip.2\Einheit
  \raise.6\Einheit\hbox to0pt{\hss.\hss}\hskip.2\Einheit
  \raise.8\Einheit\hbox to0pt{\hss.\hss}\hss}}
\def\dSSchritt{\leavevmode\raise-.4pt\hbox to0pt{%
  \hbox to0pt{\hss.\hss}\hskip.2\Einheit
  \raise-.2\Einheit\hbox to0pt{\hss.\hss}\hskip.2\Einheit
  \raise-.4\Einheit\hbox to0pt{\hss.\hss}\hskip.2\Einheit
  \raise-.6\Einheit\hbox to0pt{\hss.\hss}\hskip.2\Einheit
  \raise-.8\Einheit\hbox to0pt{\hss.\hss}\hss}}
\def\SPfad(#1,#2),#3\endSPfad{\unskip\leavevmode
  \xcoord#1 \ycoord#2 \ZeichneSPfad#3\endSPfad}
\def\ZeichneSPfad#1{\ifx#1\endSPfad\let\next\relax
  \else\let\next\ZeichneSPfad
    \ifnum#1=1
      \raise\ycoord \Einheit\hbox to0pt{\hskip\xcoord \Einheit
         \hSSchritt\hss}%
      \advance\xcoord by 1
    \else\ifnum#1=2
      \raise\ycoord \Einheit\hbox to0pt{\hskip\xcoord \Einheit
        \hbox{\hskip-2pt \vSSchritt}\hss}%
      \advance\ycoord by 1
    \else\ifnum#1=3
      \raise\ycoord \Einheit\hbox to0pt{\hskip\xcoord \Einheit
         \DSSchritt\hss}
      \advance\xcoord by 1
      \advance\ycoord by 1
    \else\ifnum#1=4
      \raise\ycoord \Einheit\hbox to0pt{\hskip\xcoord \Einheit
         \dSSchritt\hss}
      \advance\xcoord by 1
      \advance\ycoord by -1
    \fi\fi\fi\fi
  \fi\next}
\def\Koordinatenachsen(#1,#2){\unskip
 \hbox to0pt{\hskip-.5pt\vrule height#2 \Einheit width.5pt depth1 
\Einheit}%
 \hbox to0pt{\hskip-1 \Einheit \xcoord#1 \advance\xcoord by1
    \vrule height0.25pt width\xcoord \Einheit depth0.25pt\hss}}
\def\Koordinatenachsen(#1,#2)(#3,#4){\unskip
 \hbox to0pt{\hskip-.5pt \ycoord-#4 \advance\ycoord by1
    \vrule height#2 \Einheit width.5pt depth\ycoord \Einheit}%
 \hbox to0pt{\hskip-1 \Einheit \hskip#3\Einheit 
    \xcoord#1 \advance\xcoord by1 \advance\xcoord by-#3 
    \vrule height0.25pt width\xcoord \Einheit depth0.25pt\hss}}
\def\Gitter(#1,#2){\unskip \xcoord0 \ycoord0 \leavevmode
  \LOOP\ifnum\ycoord<#2
    \loop\ifnum\xcoord<#1
      \raise\ycoord \Einheit\hbox to0pt{\hskip\xcoord 
\Einheit\Punkt\hss}%
      \advance\xcoord by1
    \repeat
    \xcoord0
    \advance\ycoord by1
  \REPEAT}
\def\Gitter(#1,#2)(#3,#4){\unskip \xcoord#3 \ycoord#4 \leavevmode
  \LOOP\ifnum\ycoord<#2
    \loop\ifnum\xcoord<#1
      \raise\ycoord \Einheit\hbox to0pt{\hskip\xcoord 
\Einheit\Punkt\hss}%
      \advance\xcoord by1
    \repeat
    \xcoord#3
    \advance\ycoord by1
  \REPEAT}
\def\Label#1#2(#3,#4){\unskip \xdim#3 \Einheit \ydim#4 \Einheit
  \def\lo{\advance\xdim by-.5 \Einheit \advance\ydim by.5 \Einheit}%
  \def\llo{\advance\xdim by-.25cm \advance\ydim by.5 \Einheit}%
  \def\loo{\advance\xdim by-.5 \Einheit \advance\ydim by.25cm}%
  \def\o{\advance\ydim by.25cm}%
  \def\ro{\advance\xdim by.5 \Einheit \advance\ydim by.5 \Einheit}%
  \def\rro{\advance\xdim by.25cm \advance\ydim by.5 \Einheit}%
  \def\roo{\advance\xdim by.5 \Einheit \advance\ydim by.25cm}%
  \def\l{\advance\xdim by-.30cm}%
  \def\r{\advance\xdim by.30cm}%
  \def\lu{\advance\xdim by-.5 \Einheit \advance\ydim by-.6 \Einheit}%
  \def\llu{\advance\xdim by-.25cm \advance\ydim by-.6 \Einheit}%
  \def\luu{\advance\xdim by-.5 \Einheit \advance\ydim by-.30cm}%
  \def\u{\advance\ydim by-.30cm}%
  \def\ru{\advance\xdim by.5 \Einheit \advance\ydim by-.6 \Einheit}%
  \def\rru{\advance\xdim by.25cm \advance\ydim by-.6 \Einheit}%
  \def\ruu{\advance\xdim by.5 \Einheit \advance\ydim by-.30cm}%
  #1\raise\ydim\hbox to0pt{\hskip\xdim
     \vbox to0pt{\vss\hbox to0pt{\hss$#2$\hss}\vss}\hss}%
}
\catcode`\@=12

\usepackage{amsmath,amsthm}
\usepackage{amssymb}
\usepackage{color}
\usepackage[colorlinks=true,
linkcolor=webcyan,
filecolor=webbrown,
citecolor=red]{hyperref}
\newcommand{\textrma}[1]{\textrm{\small{#1}} }
\def\red{\textcolor{red} }
\def\blue{\textcolor{blue} }
\def\green{\textcolor{green} }
\def\v{\vert}
\def\g{\ensuremath{\mathcal{ G}_{n}}}
\def\d{\ensuremath{\mathcal{ D}_{n}}}

\def\p{\ensuremath{\mathcal P}}
\def\gs{\ensuremath{\mathcal{SG}_{n}}}
\def\ds{\ensuremath{\mathcal{SD}_{n}}}

\begin{document}
\newtheorem{lemma}{Lemma}
\newtheorem{theorem}{Theorem}
\newtheorem{prop}{Proposition}
\newtheorem{cor}{Corollary}
\begin{center}
{\Large
  A Uniformly Distributed Parameter on a Class of Lattice Paths                          \\ 
}
\vspace{10mm}
DAVID CALLAN  \\
Department of Statistics  \\
University of Wisconsin-Madison  \\
1210 W. Dayton St   \\
Madison, WI \ 53706-1693  \\
{\bf callan@stat.wisc.edu}  \\
\vspace{5mm}
\end{center}

\begin{abstract}
Let $\g$ denote the set of lattice paths from $(0,0)$ to 
$(n,n)$ with steps of the form $(i,j)$ where $i$ and $j$ are nonnegative 
integers, not both 0.  
Let $\d$ denote the set of paths in $\g$ 
with steps restricted to 
$(1,0),(0,1),(1,1)$, so-called Delannoy 
paths. Stanley has shown that $\v \g \v =2^{n-1} \v 
\d \v$ and Sulanke has given a bijective proof. 
Here we give a simple parameter on $\g$ that is uniformly 
distributed over the $2^{n-1}$ subsets of $[n-1]=\{1,2,\ldots,n\}$ and takes the 
value $[n-1]$ precisely on the Delannoy paths.	
\end{abstract}


We consider paths in the lattice plane $\mathbb{Z}^{2}$ with 
arbitrary nonnegative-integer-coordinate steps, that is, steps in 
$\mathbb{N} \times \mathbb{N} \backslash\{ (0,0)\}$, called 
\emph{general} lattice paths. A path can be specified  by the sequence 
of its steps or, depending on where the path is situated in  
$\mathbb{Z}^{2}$, either by its vertices or
by its line segments. 
Let $\g$ denote the set of general lattice paths from $(0,0)$ to 
$(n,n)$. Let $\d$ denote the set of paths in $\g$ 
with steps restricted to $(1,0),(0,1),(1,1)$, so-called Delannoy 
paths. Stanley \cite[Ex. 6.16]{ec2} shows that $\v \g \v =2^{n-1} \v 
\d \v$ and Sulanke \cite{sulanke2000} has given a bijective proof. 
Here we give a simple parameter on $\g$ that is uniformly 
distributed over the $2^{n-1}$ subsets of $[n-1]$ and takes the 
value $[n-1]$ precisely on the Delannoy paths.


To present this parameter, two notions are relevant: a path is 
\emph{balanced} if its terminal vertex lies on the line of slope 1 
through its initial vertex. A path is \emph{subdiagonal} if 
it never rises above the line of slope 1 through its initial vertex. 
Analogously for \emph{superdiagonal}. 
A \emph{subpath} of a path $\pi$ is of course a subsequence of 
consecutive steps of $\pi$. Since subpaths that do not start at the 
origin will arise, the reader should not confuse a path's 
inherent property of being subdiagonal with its placement relative to 
the diagonal line $y=x$. 

For $\pi \in \g$, consider the interior vertical lines: $x=k,\ 1 \le k 
\le n-1$. Such a line is \emph{active} for $\pi$ if it contains a 
vertex of $\pi$---an \emph{active} vertex---that (i) lies on the line $y=x$, or 
(ii) lies strictly below $y=x$ and is the initial vertex of a nonempty balanced 
subdiagonal subpath of $\pi$, or (iii) lies strictly above $y=x$ and is the 
terminal vertex of a nonempty balanced superdiagonal subpath of $\pi$. If a line 
is active for $\pi$ by virtue of (i), no other vertex on the 
line can meet the conditions of (ii) or (iii). If active by virtue of 
(ii) or (iii), then all path vertices on the line lie strictly to one 
side of $y=x$ and only the one closest to $y=x$ is active. In any 
case, an active line contains a unique active vertex.

\begin{prop}
	A path $\pi \in \g$ is Delannoy if and only if \emph{all} its interior vertical 
	lines are active.
\end{prop}
\begin{proof}
	The ``if'' part is clear. For the converse, 
suppose all lines are active for $\pi$. If $\pi$ had a line segment 
of slope $m$ with $0<m<1$, then there would be an interior vertical line 
containing no vertex of $\pi$ at all, giving an inactive line. If $\pi$ had a line 
segment $PQ$ of slope $m$ with $1<m<\infty$, then either $P$ is 
strictly below $y=x$ making the vertical line through $P$ inactive or 
$Q$ is strictly above $y=x$ likewise giving an inactive line (or 
both).

Hence all line segments in $\pi$ have slope $0,1$ or $\infty$. But a 
missing interior lattice point in 
a segment whose slope is 0 or 1 in $\pi$ would clearly give an inactive 
line. And a segment of slope $\infty$ and length $\ge 2$ in $\pi$  would 
contain an interior lattice point $P$ either (i) on $y=x$, (ii) 
above $y=x$, or (iii) below $y=x$. In case (i), the vertical line 
through $P$ is inactive. In case (ii), let $L$ denote the ray from 
$P$ heading Northeast. If $L$ contains a vertex of the path, the first 
such determines an inactive line. Otherwise, $\pi$ must cross $L$ on 
a segment $PQ$ of slope $<1$; all lines strictly between $P$ and 
$Q$ are inactive and there is at least one such. Similarly, case 
(iii) gives an inactive line. Hence all line segments in $\pi$ have slope 
$0,1$ or $\infty$ and minimum possible length; in other words, $\pi$ 
is Delannoy.
\end{proof}

The \emph{active set} for $\pi\in \g$ is $\{k\in
[n-1] \mid \,x=k$ is active for $\pi\}$.
Thus for $\pi\in \g$, Proposition 1 asserts that its active set is $[n-1]$ iff $\pi$ is Delannoy.
Our main result is 
\begin{theorem}
	The parameter ``active set'' on $\g$ is uniformly distributed over 
	all subsets of $[n-1]$.
\end{theorem}
We will first prove a similar result for subdiagonal paths. Let
$\gs$ and $\ds$ denote the set of subdiagonal paths in $\g$ and $\d$ 
respectively.
Of course, as for $\g,\ \pi\in \gs$ is Delannoy if and only if all its interior vertical 
lines are active.
\begin{theorem}
	The parameter ``active set'' restricted to $\gs$ is also
	uniformly distributed over all subsets of $[n-1]$.
\end{theorem}
To establish Theorem 2  we will define a map $f$ that takes a path 
in $\gs$ together with an active line $x=k$ and produces a path 
in $\gs$ in which $x=k$ is not active without disturbing the activity 
status of other lines: it ``deactivates'' $k$. 
The map merely deletes the 
active vertex for $k$ and adjusts the location of some of its 
successors along
the vertical line they lie on.  Furthermore, the map 
is commutative: given $k,\ell$ active for $\pi\in \gs$, you get the same 
result deactivating them in either order. The tricky part is verifying 
that $f$ is reversible. Once this is done, we have a bijective 
correspondence between $\gs$ and $\p([n-1])\times \ds$ via ``record the 
active set for $\pi\in \gs$ and then activate all of $\pi$'s inactive lines'', 
and Theorem 2 follows.

It is convenient to use the following abbreviations. For lattice 
points $X,Y,\ m_{XY}$ denotes the slope of the line $XY$ (possibly 
$\infty$). It is also convenient to assume $m_{XX}=\infty$ ($m_{XX}=1$ 
would work just as well, but not $m_{XX}<1$). Also, $L_{X}$ denotes 
the line through $X$ of slope 1 .

To define $f$, suppose given a subdiagonal path $\pi$ in $\gs$ and an 
active $k$. Locate the active vertex $P$ on $x=k$, 
its predecessor vertex $A$ 
and its successor vertex $B$ on $\pi$. 
There are two cases.

\textbf{Case} $\mathbf{m_{AP}<1}$. \quad Find the first vertex $Q$ on $\pi$ that is 
strictly above $L_{P}$. The existence 
of $Q$ is guaranteed because $m_{AP}<1$.  Lower $\pi$'s vertices $B$ 
through (= up to and including) the predecessor of $Q$ 
by $h$ units vertically where $h=$ vertical distance from $P$ down to
$B$ (possibly 0). Note that the predecessor of $Q$ may be $B$. Finally, delete $P$.

\textbf{Case} $\mathbf{m_{AP}\ge 1}$ (includes $m_{AP}=\infty$). \quad Find the first 
vertex $Q$ on $\pi$ that terminates a 
nonempty balanced subdiagonal subpath starting at $P$. Note that $Q$ 
may be $B$. Lower $\pi$'s vertices $B$ through the predecessor of $Q$ 
on $\pi$ by $h$ units vertically where now $h=$ vertical distance from 
$P$ down to 
$L_{A}$ (not to $A$). Again $h$ may $=0$ but, unlike Case 1, even if $h>0$ the set of 
vertices to be lowered may be vacuous---this occurs if $Q=B$---and 
then no vertices actually get lowered.
Finally, delete $P$.

Figure 1 gives an example of Case $m_{AP}<1$, Figure 2 gives an example of 
Case $m_{AP}\ge 1$, and Figure 3 gives the action 
of $f$ on all 6 paths in $\mathcal{SG}_{2}$ for which $x=1$ is 
active, that is, on $\mathcal{SD}_{2}$.
The active line is in red and becomes a blue inactive line in the 
image path. The unaltered vertex $A$ is evident in the image path and the 
vertex $B$ is readily recovered as the first vertex strictly to the 
right of the now inactive line. Also marked is the projection $B'$ of $B$ 
on the inactive line; $B'$ is key to reversing $f$.

To reverse $f$, we must first distinguish cases in the image path. In 
fact, we find $m_{AB'}<1$ in 
Case $m_{AP}<1$, and $m_{AB'}\ge 1$ (includes $A=B'$) in Case $m_{AP}\ge 1$. 
Then proceed as follows.


For Case $m_{AB'}<1$:
\vspace*{-3mm}
\begin{itemize}
	\item  Retrieve $Q$ as the first vertex after $A$ on the image 
           strictly above the line $L_{B'}$. 

	\item  Retrieve $h$ as the vertical distance from $L_{B'}$ down to 
           $Q$'s predecessor. Raise vertices $B$ through $Q$ by $h$ units.

	\item  Retrieve $P$ as $B'$.

\end{itemize}

For Case $m_{AB'}\ge 1$:
\vspace*{-3mm}
\begin{itemize}
	\item  Retrieve $Q$ as the first vertex strictly after $A$ that lies 
        	weakly above  $L_{A}$.

	\item  Retrieve $h$ as the vertical distance from $Q$ down to $L_{A}$. 
	       Raise vertices $B$ through the predecessor of $Q$ by $h$ units.

	\item  Retrieve $P$ at height $h$ above $L_{A}$ on the inactive line.
\end{itemize}
\Einheit=0.5cm
\[
\green{
\SPfad(-13,0),333333333333\endSPfad
\DuennPunkt(-12,0)
\DuennPunkt(-12,1)
\DuennPunkt(-11,0)
\DuennPunkt(-11,1)
\DuennPunkt(-11,2)
\DuennPunkt(-10,0)
\DuennPunkt(-10,1)
\DuennPunkt(-10,2)
\DuennPunkt(-10,3)
\DuennPunkt(-9,0)
\DuennPunkt(-9,1)
\DuennPunkt(-9,2)
\DuennPunkt(-9,3)
\DuennPunkt(-9,4)
\DuennPunkt(-8,0)
\DuennPunkt(-8,1)
\DuennPunkt(-8,2)
\DuennPunkt(-8,3)
\DuennPunkt(-8,4)
\DuennPunkt(-8,5)
\DuennPunkt(-7,0)
\DuennPunkt(-7,1)
\DuennPunkt(-7,2)
\DuennPunkt(-7,3)
\DuennPunkt(-7,4)
\DuennPunkt(-7,5)
\DuennPunkt(-7,6)
\DuennPunkt(-6,0)\DuennPunkt(-6,1)\DuennPunkt(-6,2)\DuennPunkt(-6,3)
\DuennPunkt(-6,4)\DuennPunkt(-6,5)\DuennPunkt(-6,6)\DuennPunkt(-6,7)
\DuennPunkt(-5,0)
\DuennPunkt(-5,1)
\DuennPunkt(-5,2)
\DuennPunkt(-5,3)
\DuennPunkt(-5,4)
\DuennPunkt(-5,5)
\DuennPunkt(-5,6)
\DuennPunkt(-5,7)
\DuennPunkt(-5,8)
\DuennPunkt(-4,0)
\DuennPunkt(-4,1)
\DuennPunkt(-4,2)
\DuennPunkt(-4,3)
\DuennPunkt(-4,4)
\DuennPunkt(-4,5)
\DuennPunkt(-4,6)
\DuennPunkt(-4,7)
\DuennPunkt(-4,8)
\DuennPunkt(-4,9)
\DuennPunkt(-3,0)
\DuennPunkt(-3,1)
\DuennPunkt(-3,2)
\DuennPunkt(-3,3)
\DuennPunkt(-3,4)
\DuennPunkt(-3,5)
\DuennPunkt(-3,6)
\DuennPunkt(-3,8)
\DuennPunkt(-3,9)
\DuennPunkt(-3,10)
\DuennPunkt(-2,0)
\DuennPunkt(-2,1)
\DuennPunkt(-2,2)
\DuennPunkt(-2,3)
\DuennPunkt(-2,4)
\DuennPunkt(-2,5)
\DuennPunkt(-2,6)
\DuennPunkt(-2,7)
\DuennPunkt(-2,8)
\DuennPunkt(-2,9)
\DuennPunkt(-2,10)
\DuennPunkt(-2,11)
\DuennPunkt(-1,0)
\DuennPunkt(-1,1)
\DuennPunkt(-1,2)
\DuennPunkt(-1,3)
\DuennPunkt(-1,4)
\DuennPunkt(-1,5)
\DuennPunkt(-1,6)
\DuennPunkt(-1,7)
\DuennPunkt(-1,8)
\DuennPunkt(-1,9)
\DuennPunkt(-1,10)
\DuennPunkt(1,0)
\DuennPunkt(2,0)
\DuennPunkt(3,0)
\DuennPunkt(4,0)
\DuennPunkt(5,0)
\DuennPunkt(6,0)
\DuennPunkt(7,0)
\DuennPunkt(8,0)
\DuennPunkt(9,0)
\DuennPunkt(10,0)
\DuennPunkt(11,0)
\DuennPunkt(12,0)
\DuennPunkt(13,0)
\DuennPunkt(13,1)
\DuennPunkt(13,2)
\DuennPunkt(13,3)
\DuennPunkt(13,4)
\DuennPunkt(13,5)
\DuennPunkt(13,6)
\DuennPunkt(13,7)
\DuennPunkt(13,8)
\DuennPunkt(13,9)
\DuennPunkt(13,10)
\DuennPunkt(12,0)
\DuennPunkt(12,1)
\DuennPunkt(12,2)
\DuennPunkt(12,3)
\DuennPunkt(12,4)
\DuennPunkt(12,5)
\DuennPunkt(12,6)
\DuennPunkt(12,7)
\DuennPunkt(12,8)
\DuennPunkt(12,9)
\DuennPunkt(12,10)
\DuennPunkt(12,11)
\DuennPunkt(11,0)
\DuennPunkt(11,1)
\DuennPunkt(11,2)
\DuennPunkt(11,3)
\DuennPunkt(11,4)
\DuennPunkt(11,5)
\DuennPunkt(11,6)
\DuennPunkt(11,7)
\DuennPunkt(11,8)
\DuennPunkt(11,9)
\DuennPunkt(11,10)
\DuennPunkt(10,0)
\DuennPunkt(10,1)
\DuennPunkt(10,2)
\DuennPunkt(10,3)
\DuennPunkt(10,4)
\DuennPunkt(10,5)
\DuennPunkt(10,6)
\DuennPunkt(10,7)
\DuennPunkt(10,8)
\DuennPunkt(10,9)
\DuennPunkt(9,0)
\DuennPunkt(9,1)
\DuennPunkt(9,2)
\DuennPunkt(9,3)
\DuennPunkt(9,4)
\DuennPunkt(9,5)
\DuennPunkt(9,6)
\DuennPunkt(9,7)
\DuennPunkt(9,8)
\DuennPunkt(8,0)
\DuennPunkt(8,1)
\DuennPunkt(8,2)
\DuennPunkt(8,3)
\DuennPunkt(8,4)
\DuennPunkt(8,5)
\DuennPunkt(8,6)
\DuennPunkt(8,7)
\DuennPunkt(7,0)
\DuennPunkt(7,1)
\DuennPunkt(7,2)
\DuennPunkt(7,3)
\DuennPunkt(7,4)
\DuennPunkt(7,5)
\DuennPunkt(7,6)
\DuennPunkt(6,0)
\DuennPunkt(6,1)
\DuennPunkt(6,2)
\DuennPunkt(6,3)
\DuennPunkt(6,4)
\DuennPunkt(6,5)
\DuennPunkt(5,0)
\DickPunkt(5,1)
\DuennPunkt(5,2)
\DuennPunkt(5,3)
\DuennPunkt(5,4)
\DuennPunkt(4,0)
\DuennPunkt(4,1)
\DuennPunkt(4,2)
\DuennPunkt(4,3)
\DuennPunkt(3,0)
\DuennPunkt(3,1)
\DuennPunkt(3,2)
\DuennPunkt(2,0)
\DuennPunkt(2,1)
\DuennPunkt(1,0)
\SPfad(1,0),333333333333\endSPfad
 }
\red{\Pfad(-9,0),2222\endPfad}
\SPfad(-9,1),33333333\endSPfad
\NormalPunkt(-13,0)
\NormalPunkt(-11,0)
\NormalPunkt(-9,1)
\NormalPunkt(-6,3)
\NormalPunkt(-4,4)
\NormalPunkt(-3,6)
\NormalPunkt(-3,7)
\NormalPunkt(-2,8)
\NormalPunkt(-1,10)
\NormalPunkt(-1,11)
\NormalPunkt(-1,12)
\NormalPunkt(1,0)
\NormalPunkt(3,0)
\NormalPunkt(8,1)
\NormalPunkt(10,2)
\NormalPunkt(11,4)
\NormalPunkt(11,5)
\NormalPunkt(12,6)
\NormalPunkt(13,10)
\NormalPunkt(13,11)
\NormalPunkt(13,12)
\blue{\Pfad(5,0),2222\endPfad}
\Label\o{\longrightarrow}(1,6)
\Label\o{f}(1,7)
\Label\lu{\scriptstyle{P}}(-8.9,1.1)
\Label\lu{\scriptstyle{A}}(-10.9,0.1)
\Label\ru{\scriptstyle{Q}}(-0.9,10.5)
\Label\lu{\scriptstyle{B}}(-5.9,3.1)
\Label\lu{\scriptstyle{A}}(3.1,0.1)
\Label\lu{\scriptstyle{B}}(8.1,1.1)
\Label\ru{\scriptstyle{B'}}(4.,1.1)
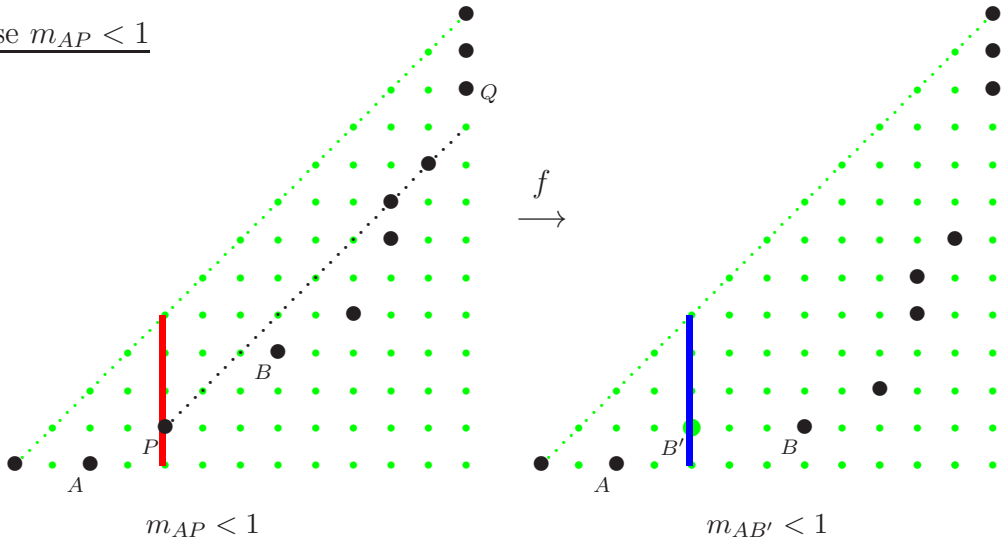
\Label\u{\textrma{vertices of path in \gs, active line in}}(-8,-2.5)
\Label\u{\textrma{red, active vertex $P$, predecessor $A$,}}(-8,-3.5)
\Label\u{\textrma{successor $B$, $Q$ first vertex after $P$ }}(-8,-4.5)
\Label\u{\textrma{strictly above $L_{P}$, predecessor of $Q$ }}(-8,-5.5)
\Label\u{\textrma{is on $L_{P}$ because $P$ is active, $h$ is }}(-8,-6.5)
\Label\u{\textrma{vertical distance from $B$ down}}(-8,-7.5)
\Label\u{\textrma{to $B$ (here 2) }}(-8,-8.5)
\Label\u{\textrma{vertices of image path, ``deactivated''}}(7,-2.5)
\Label\u{\textrma{line in blue, $B'$ projection of $B$ on }}(7,-3.5)
\Label\u{\textrma{blue line, $Q$ can be retrieved as first}}(7,-4.5)
\Label\u{\textrma{vertex strictly above $L_{B'}$, then $h$ can }}(7,-5.5)
\Label\u{\textrma{be retrieved as vertical distance}}(7,-6.5)
\Label\u{\textrma{from $L_{B'}$ down to $Q$'s predecessor}}(7,-7.5)
\Label\u{\textrma{and, lastly, $P$ is $B'$ }}(7,-8.5)
\Label\u{\textrm{\underline{Case $m_{AP}<1$}}}(-12,12)
\Label\u{\textrma{$m_{AP}< 1$}}(-8,-1.0)
\Label\u{\textrma{$m_{AB'}< 1$}}(7,-1.0)
\Label\u{\textrm{Figure 1}}(0,-10.5)
\]

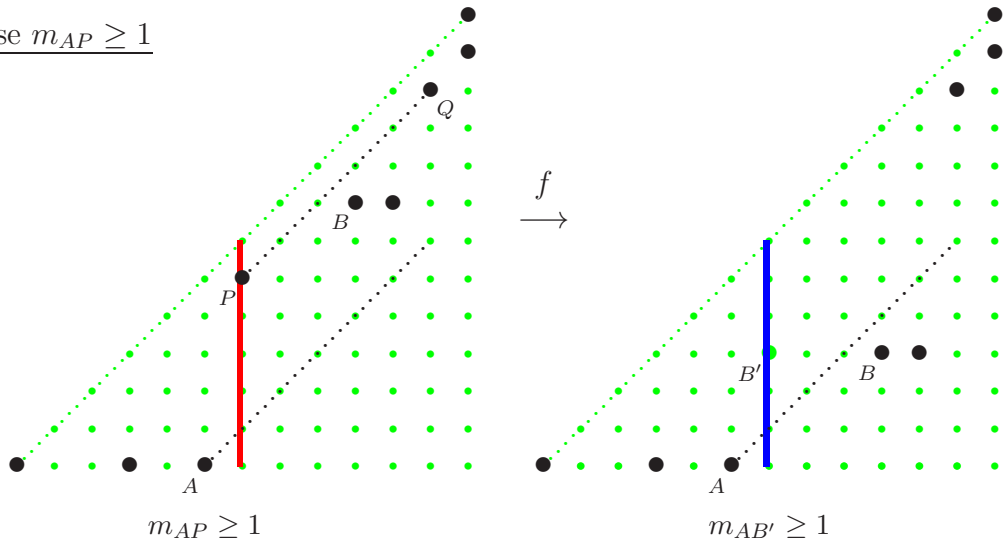
\Einheit=0.5cm
\[
\green{
\SPfad(-13,0),333333333333\endSPfad
\DuennPunkt(-12,0)
\DuennPunkt(-12,1)
\DuennPunkt(-11,0)
\DuennPunkt(-11,1)
\DuennPunkt(-11,2)
\DuennPunkt(-10,1)
\DuennPunkt(-10,2)
\DuennPunkt(-10,3)
\DuennPunkt(-9,0)
\DuennPunkt(-9,1)
\DuennPunkt(-9,2)
\DuennPunkt(-9,3)
\DuennPunkt(-9,4)
\DuennPunkt(-8,1)
\DuennPunkt(-8,2)
\DuennPunkt(-8,3)
\DuennPunkt(-8,4)
\DuennPunkt(-8,5)
\DuennPunkt(-7,0)
\DuennPunkt(-7,1)
\DuennPunkt(-7,2)
\DuennPunkt(-7,3)
\DuennPunkt(-7,4)
\DuennPunkt(-7,6)
\DuennPunkt(-6,0)\DuennPunkt(-6,1)\DuennPunkt(-6,2)\DuennPunkt(-6,3)
\DuennPunkt(-6,4)\DuennPunkt(-6,5)\DuennPunkt(-6,6)\DuennPunkt(-6,7)
\DuennPunkt(-5,0)
\DuennPunkt(-5,1)
\DuennPunkt(-5,2)
\DuennPunkt(-5,3)
\DuennPunkt(-5,4)
\DuennPunkt(-5,5)
\DuennPunkt(-5,6)
\DuennPunkt(-5,7)
\DuennPunkt(-5,8)
\DuennPunkt(-4,0)
\DuennPunkt(-4,1)
\DuennPunkt(-4,2)
\DuennPunkt(-4,3)
\DuennPunkt(-4,4)
\DuennPunkt(-4,5)
\DuennPunkt(-4,6)
\DuennPunkt(-4,8)
\DuennPunkt(-4,9)
\DuennPunkt(-3,0)
\DuennPunkt(-3,1)
\DuennPunkt(-3,2)
\DuennPunkt(-3,3)
\DuennPunkt(-3,4)
\DuennPunkt(-3,5)
\DuennPunkt(-3,6)
\DuennPunkt(-3,8)
\DuennPunkt(-3,9)
\DuennPunkt(-3,10)
\DuennPunkt(-2,0)
\DuennPunkt(-2,1)
\DuennPunkt(-2,2)
\DuennPunkt(-2,3)
\DuennPunkt(-2,4)
\DuennPunkt(-2,5)
\DuennPunkt(-2,6)
\DuennPunkt(-2,7)
\DuennPunkt(-2,8)
\DuennPunkt(-2,9)
\DuennPunkt(-2,11)
\DuennPunkt(-1,0)
\DuennPunkt(-1,1)
\DuennPunkt(-1,2)
\DuennPunkt(-1,3)
\DuennPunkt(-1,4)
\DuennPunkt(-1,5)
\DuennPunkt(-1,6)
\DuennPunkt(-1,7)
\DuennPunkt(-1,8)
\DuennPunkt(-1,9)
\DuennPunkt(-1,10)
\DuennPunkt(1,0)
\DuennPunkt(2,0)
\DuennPunkt(3,0)
\DuennPunkt(4,0)
\DuennPunkt(5,0)
\DuennPunkt(6,0)
\DuennPunkt(7,0)
\DuennPunkt(8,0)
\DuennPunkt(9,0)
\DuennPunkt(10,0)
\DuennPunkt(11,0)
\DuennPunkt(12,0)
\DuennPunkt(13,0)
\DuennPunkt(13,1)
\DuennPunkt(13,2)
\DuennPunkt(13,3)
\DuennPunkt(13,4)
\DuennPunkt(13,5)
\DuennPunkt(13,6)
\DuennPunkt(13,7)
\DuennPunkt(13,8)
\DuennPunkt(13,9)
\DuennPunkt(13,10)
\DuennPunkt(12,0)
\DuennPunkt(12,1)
\DuennPunkt(12,2)
\DuennPunkt(12,3)
\DuennPunkt(12,4)
\DuennPunkt(12,5)
\DuennPunkt(12,6)
\DuennPunkt(12,7)
\DuennPunkt(12,8)
\DuennPunkt(12,9)
\DuennPunkt(12,11)
\DuennPunkt(11,0)
\DuennPunkt(11,1)
\DuennPunkt(11,2)
\DuennPunkt(11,3)
\DuennPunkt(11,4)
\DuennPunkt(11,5)
\DuennPunkt(11,6)
\DuennPunkt(11,7)
\DuennPunkt(11,8)
\DuennPunkt(11,9)
\DuennPunkt(11,10)
\DuennPunkt(10,0)
\DuennPunkt(10,1)
\DuennPunkt(10,2)
\DuennPunkt(10,3)
\DuennPunkt(10,4)
\DuennPunkt(10,5)
\DuennPunkt(10,6)
\DuennPunkt(10,7)
\DuennPunkt(10,8)
\DuennPunkt(10,9)
\DuennPunkt(9,0)
\DuennPunkt(9,1)
\DuennPunkt(9,2)
\DuennPunkt(9,3)
\DuennPunkt(9,4)
\DuennPunkt(9,5)
\DuennPunkt(9,6)
\DuennPunkt(9,7)
\DuennPunkt(9,8)
\DuennPunkt(8,0)
\DuennPunkt(8,1)
\DuennPunkt(8,2)
\DuennPunkt(8,3)
\DuennPunkt(8,4)
\DuennPunkt(8,5)
\DuennPunkt(8,6)
\DuennPunkt(8,7)
\DuennPunkt(7,0)
\DuennPunkt(7,1)
\DuennPunkt(7,2)
\NormalPunkt(7,3)
\DuennPunkt(7,4)
\DuennPunkt(7,5)
\DuennPunkt(7,6)
\DuennPunkt(6,0)
\DuennPunkt(6,1)
\DuennPunkt(6,2)
\DuennPunkt(6,3)
\DuennPunkt(6,4)
\DuennPunkt(6,5)
\DuennPunkt(5,0)
\DuennPunkt(5,1)
\DuennPunkt(5,2)
\DuennPunkt(5,3)
\DuennPunkt(5,4)
\DuennPunkt(4,0)
\DuennPunkt(4,1)
\DuennPunkt(4,2)
\DuennPunkt(4,3)
\DuennPunkt(3,0)
\DuennPunkt(3,1)
\DuennPunkt(3,2)
\DuennPunkt(2,0)
\DuennPunkt(2,1)
\DuennPunkt(1,0)
\SPfad(1,0),333333333333\endSPfad}
\red{\Pfad(-7,0),222222\endPfad}
\SPfad(-7,5),33333\endSPfad
\NormalPunkt(-13,0)
\NormalPunkt(-10,0)
\NormalPunkt(-8,0)
\NormalPunkt(-7,5)
\NormalPunkt(-4,7)
\NormalPunkt(-3,7)
\NormalPunkt(-2,10)
\NormalPunkt(-1,11)
\NormalPunkt(-1,12)
\NormalPunkt(1,0)
\NormalPunkt(4,0)
\NormalPunkt(6,0)
\NormalPunkt(10,3)
\NormalPunkt(11,3)
\NormalPunkt(12,10)
\NormalPunkt(13,11)
\NormalPunkt(13,12)
\blue{\Pfad(7,0),222222\endPfad}
\Label\o{\longrightarrow}(1,6)
\Label\o{f}(1,7)
\Label\lu{\scriptstyle{P}}(-6.9,5.1)
\Label\lu{\scriptstyle{A}}(-7.9,0.1)
\Label\ru{\scriptstyle{Q}}(-2.1,10.1)
\Label\lu{\scriptstyle{B}}(-3.9,7.1)
\Label\lu{\scriptstyle{A}}(6.1,0.1)
\Label\lu{\scriptstyle{B}}(10.1,3.1)
\Label\ru{\scriptstyle{B'}}(6.,3.1)
\SPfad(6,0), 333333\endSPfad
\SPfad(-8,0),333333\endSPfad
\Label\u{\textrma{vertices of path in \gs, active line in}}(-8,-2.5)
\Label\u{\textrma{red, active vertex $P$, predecessor $A$,}}(-8,-3.5)
\Label\u{\textrma{successor $B$, $Q$ first vertex on $L_{P}$,}}(-8,-4.5)
\Label\u{\textrma{$h$ is vertical distance from }}(-8,-5.5)
\Label\u{\textrma{$P$ down to $L_{A}$}}(-8,-6.5)
\Label\u{\textrma{vertices of image path, ``deactivated''}}(7,-2.5)
\Label\u{\textrma{line in blue, $B'$ projection of $B$ on blue}}(7,-3.5)
\Label\u{\textrma{line, $Q$ can be retrieved as first vertex}}(7,-4.5)
\Label\u{\textrma{weakly above $L_{A}$, $h$ can be retrieved}}(7,-5.5)
\Label\u{\textrma{as vertical distance from $Q$ down}}(7,-6.5)
\Label\u{\textrma{to $L_{A}$, and $P$ is $h$ units above $L_{A}$}}(7,-7.5)
\Label\u{\textrm{\underline{Case $m_{AP}\ge 1$}}}(-12,12)
\Label\u{\textrma{$m_{AP}\ge 1$}}(-8,-1.0)
\Label\u{\textrma{$m_{AB'}\ge 1$}}(7,-1.0)
\Label\u{\textrm{Figure 2}}(0,-9.5)
\]

\vspace*{7mm}

\Einheit=0.7cm
\[
\green{
\SPfad(-11,5),33\endSPfad
\SPfad(-7,5),33\endSPfad
\SPfad(-3,5),33\endSPfad
\SPfad(1,5),33\endSPfad
\SPfad(5,5),33\endSPfad
\SPfad(9,5),33\endSPfad
\SPfad(-11,0),33\endSPfad
\SPfad(-7,0),33\endSPfad
\SPfad(-3,0),33\endSPfad
\SPfad(1,0),33\endSPfad
\SPfad(5,0),33\endSPfad
\SPfad(9,0),33\endSPfad
\DuennPunkt(-11,5)
\DuennPunkt(-11,0)
\DuennPunkt(-10,0)
\DuennPunkt(-10,1)
\DuennPunkt(-10,5)
\DuennPunkt(-10,6)
\DuennPunkt(-9,0)
\DuennPunkt(-9,1)
\DuennPunkt(-9,2)
\DuennPunkt(-9,5)
\DuennPunkt(-9,6)
\DuennPunkt(-9,7)
\DuennPunkt(-7,0)
\DuennPunkt(-7,5)
\DuennPunkt(-6,0)\DuennPunkt(-6,1)
\DuennPunkt(-6,5)\DuennPunkt(-6,6)
\DuennPunkt(-5,0)
\DuennPunkt(-5,1)
\DuennPunkt(-5,2)
\DuennPunkt(-5,5)
\DuennPunkt(-5,6)
\DuennPunkt(-5,7)
\DuennPunkt(-3,0)
\DuennPunkt(-3,5)
\DuennPunkt(-2,0)
\DuennPunkt(-2,1)
\DuennPunkt(-2,5)
\DuennPunkt(-2,6)
\DuennPunkt(-1,0)
\DuennPunkt(-1,1)
\DuennPunkt(-1,2)
\DuennPunkt(-1,5)
\DuennPunkt(-1,6)
\DuennPunkt(-1,7)
\DuennPunkt(1,0)
\DuennPunkt(1,5)
\DuennPunkt(2,0)
\DuennPunkt(2,1)
\DuennPunkt(2,5)
\DuennPunkt(2,6)
\DuennPunkt(3,0)
\DuennPunkt(3,1)
\DuennPunkt(3,2)
\DuennPunkt(3,5)
\DuennPunkt(3,6)
\DuennPunkt(3,7)
\DuennPunkt(5,0)
\DuennPunkt(5,5)
\DuennPunkt(6,0)
\DuennPunkt(6,1)
\DuennPunkt(6,5)
\DuennPunkt(6,6)
\DuennPunkt(11,0)
\DuennPunkt(11,1)
\DuennPunkt(11,2)
\DuennPunkt(11,5)
\DuennPunkt(11,6)
\DuennPunkt(11,7)
\DuennPunkt(10,0)
\DuennPunkt(10,1)
\DuennPunkt(10,5)
\DuennPunkt(10,6)
\DuennPunkt(9,0)
\DuennPunkt(9,5)
\DuennPunkt(7,0)
\DuennPunkt(7,1)
\DuennPunkt(7,2)
\DuennPunkt(7,5)
\DuennPunkt(7,6)
\DuennPunkt(7,7)
 }
\red{\Pfad(-10,5),22\endPfad}
\red{\Pfad(-6,5),22\endPfad}
\red{\Pfad(-2,5),22\endPfad}
\red{\Pfad(2,5),22\endPfad}
\red{\Pfad(6,5),22\endPfad}
\red{\Pfad(10,5),22\endPfad}
\NormalPunkt(-11,5)
\NormalPunkt(-10,6)
\NormalPunkt(-9,7)
\NormalPunkt(-7,5)
\NormalPunkt(-6,6)
\NormalPunkt(-5,6)
\NormalPunkt(-5,7)
\NormalPunkt(-3,5)
\NormalPunkt(-2,5)
\NormalPunkt(-2,6)
\NormalPunkt(-1,7)
\NormalPunkt(1,5)
\NormalPunkt(2,5)
\NormalPunkt(2,6)
\NormalPunkt(3,6)
\NormalPunkt(3,7)
\NormalPunkt(5,5)
\NormalPunkt(6,5)
\NormalPunkt(7,6)
\NormalPunkt(7,7)
\NormalPunkt(9,5)
\NormalPunkt(10,5)
\NormalPunkt(11,5)
\NormalPunkt(11,6)
\NormalPunkt(11,7)
\Pfad(-11,5),33\endPfad
\Pfad(-7,5),312\endPfad
\Pfad(-3,5),123\endPfad
\Pfad(1,5),1212\endPfad
\Pfad(5,5),132\endPfad
\Pfad(9,5),1122\endPfad
\blue{\Pfad(-10,0),22\endPfad}
\blue{\Pfad(-6,0),22\endPfad}
\blue{\Pfad(-2,0),22\endPfad}
\blue{\Pfad(2,0),22\endPfad}
\blue{\Pfad(6,0),22\endPfad}
\blue{\Pfad(10,0),22\endPfad}
\Label\o{\downarrow}(-10,3)
\Label\o{\downarrow}(-6,3)
\Label\o{\downarrow}(-2,3)
\Label\o{\downarrow}(2,3)
\Label\o{\downarrow}(6,3)
\Label\o{\downarrow}(10,3)
\Label\lo{\scriptstyle{A}}(-10.9,4.8)
\Label\lo{\scriptstyle{A}}(-6.9,4.8)
\Label\lo{\scriptstyle{A}}(-1.9,4.8)
\Label\lo{\scriptstyle{A}}(2.1,4.8)
\Label\lo{\scriptstyle{A}}(5.1,4.8)
\Label\lo{\scriptstyle{A}}(9.1,4.8)
\Label\lo{\scriptstyle{A}}(-10.9,-0.2)
\Label\lo{\scriptstyle{A}}(-6.9,-0.2)
\Label\lo{\scriptstyle{A}}(-1.9,-0.2)
\Label\lo{\scriptstyle{A}}(2.1,-0.2)
\Label\lo{\scriptstyle{A}}(5.1,-0.2)
\Label\lo{\scriptstyle{A}}(9.1,-0.2)
\Label\lo{\scriptstyle{P}}(-9.9,5.8)
\Label\lo{\scriptstyle{P}}(-5.9,5.8)
\Label\lo{\scriptstyle{P}}(-1.9,5.8)
\Label\lo{\scriptstyle{P}}(2.1,5.8)
\Label\lo{\scriptstyle{P}}(6.8,4.4)
\Label\lo{\scriptstyle{P}}(10.8,4.8)
\Label\ru{\scriptstyle{B}}(10.9,5.7)
\Label\ru{\scriptstyle{B}}(6.9,7.7)
\Label\ru{\scriptstyle{B}}(2.9,6.7)
\Label\ru{\scriptstyle{B}}(-1.1,7.7)
\Label\ru{\scriptstyle{B}}(-5.1,6.7)
\Label\ru{\scriptstyle{B}}(-9.1,7.7)
\Label\ru{\scriptstyle{B}}(10.9,0.7)
\Label\ru{\scriptstyle{B}}(6.9,0.7)
\Label\ru{\scriptstyle{B}}(2.9,0.7)
\Label\ru{\scriptstyle{B}}(-1.1,2.7)
\Label\ru{\scriptstyle{B}}(-5.1,1.7)
\Label\ru{\scriptstyle{B}}(-9.1,2.7)
\NormalPunkt(-11,0)
\NormalPunkt(-9,2)
\NormalPunkt(-7,0)
\NormalPunkt(-5,1)
\NormalPunkt(-5,2)
\NormalPunkt(-3,0)
\NormalPunkt(-2,0)
\NormalPunkt(-1,2)
\NormalPunkt(1,0)
\NormalPunkt(2,0)
\NormalPunkt(3,0)
\NormalPunkt(3,2)
\NormalPunkt(5,0)
\NormalPunkt(7,0)
\NormalPunkt(7,2)
\NormalPunkt(9,0)
\NormalPunkt(11,0)
\NormalPunkt(11,1)
\NormalPunkt(11,2)
\Label\u{\textrma{action of $f$ on  $\mathcal{SD}_{2}$}}(0,9)
\Label\u{\textrm{Figure 3}}(0,-1.5)
\]

\newpage

Having defined the bijection $f:\{(\pi,k): \pi \in \gs,\ x=k$ active 
for $\pi\} \longrightarrow \{(\pi,k): \pi \in \gs,\ x=k$ inactive 
for $\pi\}$ to establish Theorem 2, we need only extend $f$ from its 
definition on $\gs$ 
to $\tilde{f}$ defined on $\g$: Theorem 1 then follows in the same way as did 
Theorem 2.  

To do so, given $\pi \in \g$ and an active $k$, find the active vertex 
$P$ on $x=k$, its predecessor vertex $A$ and successor vertex $B$. 
We divide the possibilities into five cases.

\noindent \textbf{Case 1}\quad  \textit{$A,P,B$ all lie weakly below 
$y=x$}. If $m_{AP}<1$, $\tilde{f}$ coincides with $f$. If $m_{AP}\ge 
1$, modify $Q$ in the definition of $f$: take $Q$ as the first vertex 
strictly after $P$ that lies weakly above $L_{P}$ rather 
than ``that terminates a 
nonempty balanced subdiagonal subpath starting at $P$''. Then $\tilde{f}$ is 
defined as $f$ was. This modification is necessary because if $P$ lies 
on $y=x$, there need not be any balanced subdiagonal path $PQ$. To 
recapture the original path, $h$ is now recaptured as the minimum of
the vertical distance from $Q$ down to $L_{A}$ and the vertical distance 
from $y=x$ down to $L_{A}$.

\noindent \textbf{Case 2}\quad  \textit{$A,P,B$ all lie weakly above 
$y=x$}. Rotate everything 180$^{\circ}$ so that 
Case 1 applies, apply $\tilde{f}$, and rotate back.

\noindent \textbf{Case 3}\quad  \textit{$A$ strictly above $y=x$, $P$ 
strictly below $y=x$}. Here $m_{AP}<1$ and $\tilde{f}$ coincides with 
$f$. Note that $B$ is also strictly below $y=x$ or else $P$ would not 
be active.

\noindent \textbf{Case 4}\quad  \textit{$P$ strictly above $y=x$, $B$ 
strictly below $y=x$}. Here $A$, like $P$, is strictly above $y=x$ 
for the same reason as in Case 3. Rotate 180$^{\circ}$, apply $f$,  and 
rotate back.

\noindent \textbf{Case 5}\quad  \textit{$P$ on $y=x$, $A,B$ on
strictly opposite sides of $y=x$}. Here $\tilde{f}$ is simply 
``delete $P$''.

These five cases are exhaustive and mutually exclusive save for one 
slight overlap: if $A,P,B$ all lie on $y=x$, then Cases 1 and 2 both 
apply, but both give the same result, namely, delete $P$.

It is evident from the reversibility of $f$ that, if we know which 
case an image path arose from, we can recapture the original path. So 
we need to find distinguishing features in the image paths in the five 
cases, and to verify that every pair $(\pi,k)$ with $\pi \in \g$ and $k$ 
inactive for $\pi$ falls in one of the image cases. 

First, we can recover $A,B$ in all cases as the last vertex preceding 
the lattice point $(k,k)$ and the first vertex following $(k,k)$ 
respectively, where lattice points are ordered primarily by 
$x$-coordinate and secondarily by $y$-coordinate. The following Table 
now gives distinguishing features in the five cases.

\begin{tabular}{|c|c|c|}
	\hline
	Case & domain path $\pi$,\ \ $x=k$ active & 
	image path $\tilde{f}(\pi)$,\ \ $x=k$ inactive \\
	\hline\hline
	1 & $A,P,B$ all weakly below $y=x$ & $A,B$ weakly below $y=x$  \\
	\hline
	2 & $A,P,B$ all weakly above $y=x$ & $A,B$ weakly above $y=x$  \\
	\hline
	3 & $A$ strictly above $y=x$ and & $A$ strictly above $y=x$ and  \\
	 & $P,B$ strictly below $y=x$ & $B$ strictly below $y=k$  \\
	\hline
	4 & $A,P$ strictly above $y=x$ and & $A$ strictly above $y=k$ and  \\
	  & $B$ strictly below $y=x$  & $B$ strictly below $y=x$  \\
	\hline  
	5 & $P$ on $y=x$ and $A,B$ on strictly  & $A,B$ on strictly opposite 
	sides of $y=x$ and \\
	  &  opposite sides of $y=x$ & $A,B$ on weakly opposite 
	sides of $y=k$ \\
	\hline
\end{tabular}
\begin{center}
	How to determine the case $\tilde{f}(\pi)$ came from
\end{center}

\vspace*{3mm}

The image path cases are mutually exclusive save for the overlap in 
cases 1 and 2 when $A$ and $B$ both lie on $y=x$. Let us confirm they 
are exhaustive. If $A$ and $B$ lie weakly on the same side of $y=x$, 
then case 1 or 2 applies. Otherwise, $A$ and $B$ lie on strictly 
opposite sides of $y=x$. Now, if $A$ is strictly below $y=x$, part of 
case 5 applies (with $A$ strictly below $y=x$, $B$ strictly above 
$y=x$; this forces $A,B$ to lie weakly---in fact, strictly---on opposite 
sides of $y=k$). 

This leaves the case $A$ strictly above $y=x$ (for if $A$ 
were on $y=x$, then $A$ and $B$ would lie weakly on the same side of 
$y=x$) and $B$ strictly below $y=x$. If $A,B$ both lie strictly below 
$y=k$, Case 3 applies. If $A,B$ both lie strictly above 
$y=k$, Case 4 applies. The remaining case---$A,B$ weakly on opposite 
sides of $y=k$---is the other part of Case 5.

Figure 4 gives an example of $\tilde{f}$ (Case 1).

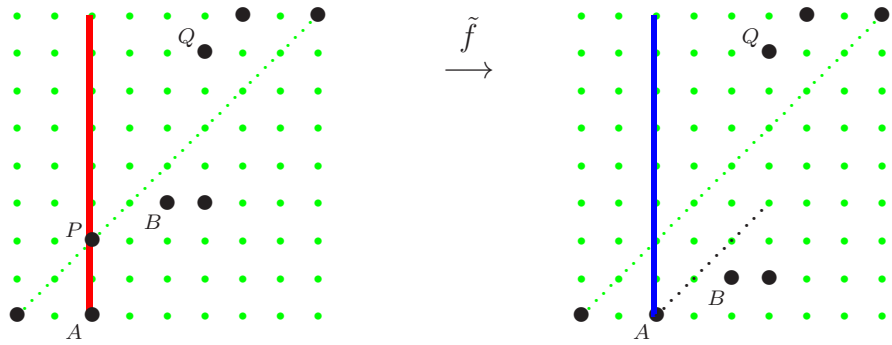
\Einheit=0.5cm
\[
\green{
\SPfad(-11,0),33333333\endSPfad
\DuennPunkt(-11,0)
\DuennPunkt(-11,1)\DuennPunkt(-11,2)\DuennPunkt(-11,3)\DuennPunkt(-11,4)
\DuennPunkt(-11,5)\DuennPunkt(-11,6)\DuennPunkt(-11,7)\DuennPunkt(-11,8)
\DuennPunkt(-10,0)\DuennPunkt(-10,1)\DuennPunkt(-10,2)\DuennPunkt(-10,3)
\DuennPunkt(-10,4)\DuennPunkt(-10,5)\DuennPunkt(-10,6)\DuennPunkt(-10,7)
\DuennPunkt(-10,8)
\DuennPunkt(-9,0)\DuennPunkt(-9,1)\DuennPunkt(-9,2)\DuennPunkt(-9,3)
\DuennPunkt(-9,4)\DuennPunkt(-9,5)\DuennPunkt(-9,6)\DuennPunkt(-9,7)
\DuennPunkt(-9,8)
\DuennPunkt(-8,0)\DuennPunkt(-8,1)\DuennPunkt(-8,2)\DuennPunkt(-8,3)
\DuennPunkt(-8,4)\DuennPunkt(-8,5)\DuennPunkt(-8,6)\DuennPunkt(-8,7)
\DuennPunkt(-8,8)
\DuennPunkt(-7,0)\DuennPunkt(-7,1)\DuennPunkt(-7,2)\DuennPunkt(-7,3)
\DuennPunkt(-7,4)\DuennPunkt(-7,5)\DuennPunkt(-7,6)\DuennPunkt(-7,7)
\DuennPunkt(-7,8)
\DuennPunkt(-6,0)\DuennPunkt(-6,1)\DuennPunkt(-6,2)\DuennPunkt(-6,3)
\DuennPunkt(-6,4)\DuennPunkt(-6,5)\DuennPunkt(-6,6)\DuennPunkt(-6,7)
\DuennPunkt(-6,8)
\DuennPunkt(-5,0)
\DuennPunkt(-5,1)
\DuennPunkt(-5,2)
\DuennPunkt(-5,3)
\DuennPunkt(-5,4)
\DuennPunkt(-5,5)
\DuennPunkt(-5,6)
\DuennPunkt(-5,7)
\DuennPunkt(-5,8)
\DuennPunkt(-4,0)
\DuennPunkt(-4,1)
\DuennPunkt(-4,2)
\DuennPunkt(-4,3)
\DuennPunkt(-4,4)
\DuennPunkt(-4,5)
\DuennPunkt(-4,6)
\DuennPunkt(-4,7)
\DuennPunkt(-4,8)
\DuennPunkt(-3,0)
\DuennPunkt(-3,1)
\DuennPunkt(-3,2)
\DuennPunkt(-3,3)
\DuennPunkt(-3,4)
\DuennPunkt(-3,5)
\DuennPunkt(-3,6)
\DuennPunkt(-3,7)
\DuennPunkt(-3,8)
\DuennPunkt(12,0)
\DuennPunkt(12,1)
\DuennPunkt(12,2)
\DuennPunkt(12,3)
\DuennPunkt(12,4)
\DuennPunkt(12,5)
\DuennPunkt(12,6)
\DuennPunkt(12,7)
\DuennPunkt(12,8)
\DuennPunkt(11,0)
\DuennPunkt(11,1)
\DuennPunkt(11,2)
\DuennPunkt(11,3)
\DuennPunkt(11,4)
\DuennPunkt(11,5)
\DuennPunkt(11,6)
\DuennPunkt(11,7)
\DuennPunkt(11,8)
\DuennPunkt(10,0)
\DuennPunkt(10,1)
\DuennPunkt(10,2)
\DuennPunkt(10,3)
\DuennPunkt(10,4)
\DuennPunkt(10,5)
\DuennPunkt(10,6)
\DuennPunkt(10,7)
\DuennPunkt(10,8)
\DuennPunkt(9,0)
\DuennPunkt(9,1)
\DuennPunkt(9,2)
\DuennPunkt(9,3)
\DuennPunkt(9,4)
\DuennPunkt(9,5)
\DuennPunkt(9,6)
\DuennPunkt(9,7)
\DuennPunkt(9,8)
\DuennPunkt(8,0)
\DuennPunkt(8,1)
\DuennPunkt(8,2)
\DuennPunkt(8,3)
\DuennPunkt(8,4)
\DuennPunkt(8,5)
\DuennPunkt(8,6)
\DuennPunkt(8,7)
\DuennPunkt(8,8)
\DuennPunkt(7,0)
\DuennPunkt(7,1)
\DuennPunkt(7,2)
\DuennPunkt(7,3)
\DuennPunkt(7,4)
\DuennPunkt(7,5)
\DuennPunkt(7,6)\DuennPunkt(7,7)\DuennPunkt(7,8)
\DuennPunkt(6,0)
\DuennPunkt(6,1)
\DuennPunkt(6,2)
\DuennPunkt(6,3)
\DuennPunkt(6,4)
\DuennPunkt(6,5)\DuennPunkt(6,6)\DuennPunkt(6,7)\DuennPunkt(6,8)
\DuennPunkt(5,0)
\DuennPunkt(5,1)
\DuennPunkt(5,2)
\DuennPunkt(5,3)\DuennPunkt(5,4)\DuennPunkt(5,5)\DuennPunkt(5,6)
\DuennPunkt(5,7)\DuennPunkt(5,8)
\DuennPunkt(4,0)\DuennPunkt(4,1)\DuennPunkt(4,2)\DuennPunkt(4,3)
\DuennPunkt(4,4)\DuennPunkt(4,5)\DuennPunkt(4,6)\DuennPunkt(4,7)
\DuennPunkt(4,8)
\SPfad(4,0),33333333\endSPfad}
\red{\Pfad(-9,0),22222222\endPfad}
\NormalPunkt(-11,0)
\NormalPunkt(-9,0)
\NormalPunkt(-9,2)
\NormalPunkt(-7,3)
\NormalPunkt(-6,3)
\NormalPunkt(-6,7)
\NormalPunkt(-5,8)
\NormalPunkt(-3,8)
\NormalPunkt(4,0)
\NormalPunkt(6,0)
\NormalPunkt(8,1)
\NormalPunkt(9,1)
\NormalPunkt(9,7)
\NormalPunkt(10,8)
\NormalPunkt(12,8)
\blue{\Pfad(6,0),22222222\endPfad}
\Label\o{\longrightarrow}(1,6)
\Label\o{\tilde{f}}(1,7)
\Label\lu{\scriptstyle{P}}(-9.0,2.9)
\Label\lu{\scriptstyle{A}}(-9.0,0.2)
\Label\ru{\scriptstyle{Q}}(-7.0,8.0)
\Label\lu{\scriptstyle{B}}(-6.9,3.1)
\Label\lu{\scriptstyle{A}}(6.1,0.2)
\Label\lu{\scriptstyle{B}}(8.1,1.1)
\Label\ru{\scriptstyle{Q}}(8.0,8.0)
\SPfad(6,0), 333\endSPfad
\Label\u{\textrma{$m_{AP}\ge 1$ and so $h$ is }}(-8,-2.5)
\Label\u{\textrma{vertical distance from $P$ }}(-8,-3.5)
\Label\u{\textrma{down to $L_{A}$ (here 2)}}(-8,-4.5)
\Label\u{\textrma{$Q$ can be retrieved as before,}}(8,-2.5)
\Label\u{\textrma{$h$ can be retrieved as minimum}}(8,-3.5)
\Label\u{\textrma{of vertical distance from $Q$ down to}}(8,-4.5)
\Label\u{\textrma{$L_{A}$ (here 4) and vertical distance }}(8,-5.5)
\Label\u{\textrma{from $y=x$ down to $L_{A}$ (here 2)}}(8,-6.5)
\Label\u{\textrm{Figure 4}}(0,-8.5)
\]

\vspace*{7mm}
This completes the proof of Theorem 1.

\end{document}